\documentclass[11pt]{amsart}
\usepackage{amsmath}
\usepackage{amsthm}
\usepackage{amsopn}
\usepackage{amssymb}
\usepackage{fullpage}
\usepackage{enumerate}
\numberwithin{equation}{section}
\usepackage{float}
\usepackage{graphicx}
\usepackage{caption}
\usepackage{subcaption}
\usepackage{titlesec}
\usepackage{textcomp}
\usepackage{fancyhdr}
\newtheorem{theorem}{Theorem}[section]
\newtheorem{proposition}{Proposition}[section]
\newtheorem{definition}{Definition}[section]
\theoremstyle{definition}
\newtheorem{remark}{Remark}[section]
\DeclareMathOperator{\Capac}{cap}
\DeclareMathOperator{\arctanh}{arctanh}
\DeclareMathOperator{\Ar}{A}
\DeclareMathOperator{\T}{\mathbb{T}}
\DeclareMathOperator{\caph}{caph}
\DeclareMathOperator{\D}{\mathbb{D}}
\DeclareMathOperator{\RE}{Re}
\DeclareMathOperator{\dist}{dist}
\DeclareMathOperator{\len}{length}
\DeclareMathOperator{\IM}{Im}
\titleformat{\section}
  {\normalfont\scshape}{\thesection}{1em}{\centering}
\titleformat{\subsection}[runin]
  {\bfseries}{\thesubsection}{1em}{}

\begin{document}
\title{COMPACT SETS IN PETALS AND THEIR BACKWARD ORBITS UNDER SEMIGROUPS OF HOLOMORPHIC FUNCTIONS}

\author{Maria Kourou}  
\author{Konstantinos Zarvalis}

\fancyhf{}
\renewcommand{\headrulewidth}{0.4pt}
\fancyhead[RO,LE]{\small \thepage}
\fancyhead[CE]{\footnotesize M. Kourou and K. Zarvalis}
\fancyhead[CO]{\footnotesize Compact Sets in Petals and their backward orbits} 
\fancyfoot[L,R,C]{}


\date{}

\begin{abstract}
Let $(\phi_t)_{t \geq 0}$ be a semigroup of holomorphic functions in the unit disk $\mathbb{D}$ and $K$ a compact subset of $\mathbb{D}$.
We investigate the conditions under which the backward orbit of $K$ under the semigroup exists. Subsequently, the geometric characteristics, as well as, potential theoretic quantities for the backward orbit of $K$ are examined. More specifically, results are obtained concerning the asymptotic behavior of its hyperbolic area and diameter, the harmonic measure and the capacity of the condenser that $K$ forms with the unit disk. 
\end{abstract}

\maketitle

\section{ \bf Introduction}

One-parameter semigroups of holomorphic self-maps of the unit disk $\D$ have been extensively examined in recent years. Their theory was introduced by Berkson and Porta in \cite{berksonporta} and later expanded in several works such as \cite{abate, BetsDescript, petals, cdmg, AnalyticFlows, shoikhetelin, goryainov, jacobzonsurvey} and in particular, in the recent monograph \cite{Booksem}. 
A one-parameter semigroup is a family $(\phi_t)_{t \geq 0}$ of holomorphic functions in $\D$, where 
\begin{enumerate}[(i)]
\item $\phi_0$ is the identity map;
\item $\phi_{t+s}(z) = \phi_t \left( \phi_s (z) \right)$, for every $t,s \geq 0$ and $z \in \D$;
\item $\phi_t(z) \xrightarrow{t \to 0^{+}} z$, uniformly on compacta in $\D$.
\end{enumerate}

One of the most important properties of one-parameter semigroups is the direct aspect of the continuous Denjoy-Wolff theorem. There exists a unique fixed point $\tau \in \overline{\D}$ such that 
\begin{equation}\label{limtraj}
\lim_{t \to +\infty} \phi_t(z)= \tau,
\end{equation}
for every point $z \in \D$.  This point $\tau$ is called the Denjoy-Wolff point of the semigroup; see \cite[Theorem 1.4.17]{abate}.
If $\tau \in\D$ and $\phi_t$ is not an elliptic automorphism of $\D$ for any $t \geq 0$, then $(\phi_t)$ is called an \textit{elliptic} semigroup. In the case where $\tau $ is a boundary point of $\D$, we observe the angular derivative $\phi_t^{\prime}(\tau)$. If $\phi_t^{\prime}(\tau)<1$, the semigroup $(\phi_t)$ is called \textit{hyperbolic}, whereas, if $\phi_t^{\prime}(\tau)=1$, the semigroup $(\phi_t)$ is called \textit{parabolic}. If $\phi_{t_0}$ is a hyperbolic (respectively parabolic) automorphism of $\D$ for some $t_0 \geq 0$, then $(\phi_t)$ is called a \textit{hyperbolic (respectively parabolic) group}. 

In the present work, we are interested in non-elliptic semigroups with boundary Denjoy-Wolff point $\tau$.
The curve $ \gamma_z: [0, + \infty)  \to \D$ with $ \gamma_z(t)= \phi_t(z)$ is called the \textit{trajectory} of a point $z \in \D$ and according to \eqref{limtraj},  $\gamma_z(t) \xrightarrow{t \to + \infty} \tau$,
for every $z \in \D$. 
So, for every point in $\D$, its trajectory approaches a fixed point in $\overline{\D}$.

In \cite{bkkp2} and \cite{kourou3}, the trajectory of a compact subset of $\D$ under a non-elliptic semigroup was examined and results on its asymptotic behavior were extracted by means of several potential theoretic and geometric quantities.
The main goal of the present work is to generalize those results through the investigation of the backward trajectory of a compact subset of $\D$ under a non-elliptic semigroup. 

Following the notation of \cite[Chapter 13]{Booksem}, the \textit{backward orbit} of a semigroup $(\phi_t)$ at a point $z \in \D$ is a continuous curve $\gamma_z: [0,+ \infty) \to \D$ that satisfies $\phi_s( \gamma_z(t))= \gamma_z(t-s)$, for every $t \in [0, + \infty)$ and every $s \in [0, t]$. 
A backward orbit is said to be \textit{regular} if $$\limsup_{t \to +\infty} d_{\D} (\gamma_z(t), \gamma_z(t+1)) < + \infty, $$
where $d_{\D}$ denotes the hyperbolic distance in $\D$. We state at this point some definitions related to backward orbits in order to meet the conditions under which the backward orbit of a compact set is defined. 

A continuous curve $\gamma : (a, + \infty)$, where $a \in [-\infty, 0)$, is called a \textit{maximal invariant curve} for $(\phi_t)$ if $$\phi_s (\gamma(t))= \gamma(t+s), \quad \forall s\geq 0, t \in (a, +\infty),  $$
$\gamma(t)\xrightarrow{t \to + \infty} \tau$ and there exists $p \in \partial \D$ such that $\gamma(t) \xrightarrow{ t \to a^{+}} p$. The point $p $ is called the \textit{starting point} of $\gamma$. From \cite[Prop. 13.3.5]{Booksem}, for every $z \in \D$, there exists a unique maximal invariant curve $\gamma_z : (a_z, +\infty)$, such that $\gamma_z (0) =z$.
The \textit{backward invariant set} $\mathcal{W}$ of $(\phi_t)$ is the set of all points in $\D$ for which $a_z = -\infty$ and it is defined as $$\mathcal{W} :=  \bigcap_{t \geq 0} \phi_t(\D) . $$
A \textit{petal} $\Delta $ of $(\phi_t)$ is a non-empty simply connected component of the interior of $\mathcal{W}$ that satisfies the following properties:
\begin{enumerate}[(i)]
   \item $\phi_t(\Delta)=\Delta $ for all $t \geq 0$ and $(\phi_{t |\Delta})$ is a group of automorphisms of $\Delta$
    \item $\tau \in \partial \Delta$
    \item there exists a boundary point $\sigma \in \partial\D \cap \partial\Delta $ (possibly $\sigma =\tau$) such that for every $z \in \Delta $ the curve $$[0,+ \infty) \ni t  \mapsto \phi_{t}^{-1}(z)$$
    is a backward regular orbit for $(\phi_t)$ that converges to $\sigma$. In addition, $\sigma $ is a boundary regular fixed point and in the case where $\sigma =\tau$, the semigroup is parabolic. 
\end{enumerate}
Moreover, as we can see in $(\text{iii})$, for any point $z \in \Delta $, we can denote, for the sake of simplicity, $\phi_{-t}(z): = \phi_{t}^{-1}(z) $ and as a result, $\phi_t(z)$ is defined for all $t \in\mathbb{R}$. More information concerning backward orbits and the characterization of petals follows in Subsection \ref{semigroups}. 

Therefore, if we suppose that $K$ is a compact subset of a petal $\Delta $ of $(\phi_t)$, then its backward trajectory is $$\gamma_K(t): = \bigcup_{z \in K} \gamma_z(t) = \phi_t^{-1}(K)=\phi_{-t}(K), \quad t \geq 0. $$
The backward orbit of every $z \in K$ is regular and so, $\gamma_K(t)$ is also regular.

With a trivial re-parametrization, we can denote the backward trajectory of $K$ by $\phi_t(K)$, $t\leq 0$. 
As $t$ decreases and tends to $- \infty$, the compact set $\phi_t (K)$ approaches the unit circle and shrinks to a boundary fixed point. 
Our purpose is to determine how the geometric and potential theoretic characteristics of $\phi_t (K)$ behave during this approach. 

We initiate our observations with the harmonic measure in the unit disk. Suppose that $(\phi_t)$ is a non-elliptic semigroup with Denjoy-Wolff point $\tau$. Further suppose that $\Delta$ is a petal of $(\phi_t)$ and $K$ is a non-polar compact subset of $\Delta$; the reader may refer to Subsection \ref{diameter} for polar sets.
The harmonic measure $\omega (\phi_t (z), \partial \phi_t  (K), \D \setminus \phi_t (K))$ is the Perron solution to the Dirichlet problem on $\D \setminus \phi_t(K)$ with given boundary values $1$ on the boundary of $\phi_t (K)$ and $0$ on the unit circle $\T$. 
Considering the harmonic measure $\omega (\phi_t (z), \partial \phi_t  (K), \D \setminus \phi_t (K))$ as a function of $t \leq 0$, we obtain the following monotonicity result. 

\begin{theorem}\label{monotharm}
Let $(\phi_t)$ be a non-elliptic semigroup of holomorphic functions in $\D$. Suppose $\Delta $ is a petal of $(\phi_t)$. Let $K$ be a compact non-polar subset of $\Delta$. 

The harmonic measure $\omega (\phi_t (z), \partial \phi_t  (K), \D \setminus \phi_t (K))$ is an increasing function of $t \in (-\infty,0]$, for every $z \in \Delta \setminus K$.
\end{theorem}

Therefore, the limit of the harmonic measure $\omega (\phi_t (z), \partial \phi_t  (K), \D \setminus \phi_t (K))$, as $t \to -\infty$, exists. This way, we can get information on the asymptotic behavior of $\phi_t (K)$. 

\begin{theorem}\label{th2}
Let $(\phi_t)$ be a non-elliptic semigroup of holomorphic functions in $\D$. Suppose $\Delta $ is a petal of $(\phi_t)$. Let $K$ be a compact non-polar subset of $\Delta$. 
Then $$\lim_{t \to - \infty} \omega (\phi_t (z), \partial \phi_t  (K), \D \setminus \phi_t (K)) = \omega (z,\partial K, \Delta\setminus K), \quad z \in \Delta \setminus K.  $$
\end{theorem}

Furthermore, we examine the change in the size of $\phi_t (K)$, as $t$ decreases and approaches $- \infty$. A natural way to do so is by observing the hyperbolic geometric characteristics of $\phi_t(K)$. 
Before we proceed to these characteristics, we need to examine the asymptotic behavior of the hyperbolic metric, as $t \to -\infty$. 
 \begin{theorem}\label{limdist}
     Let $(\phi_t)$ be a non-elliptic semigroup in $\mathbb{D}$. Suppose $\Delta$ is a petal of $(\phi_t)$. Then,
     $$\lim\limits_{t\to-\infty}\lambda_\mathbb{D}(\phi_t(z))=\lambda_\Delta(z), \quad z\in\Delta$$
     uniformly on compacta. Moreover,
     $$\lim\limits_{t\to-\infty}d_\mathbb{D}(\phi_t(z),\phi_t(w))=d_\Delta(z,w),$$
     for all $z,w\in\Delta$.
    \end{theorem}

At this point, we state a monotonicity property of the hyperbolic $n$-th diameter; see Subsection \ref{diameter}. 
\begin{theorem}\label{monothypdiam}
Suppose $(\phi_t)$ is a non-elliptic semigroup of holomorphic functions in $\D$, $\Delta$ is a petal of $(\phi_t)$ and $K \subset \Delta$ compact.
The hyperbolic $n$-th diameter $ d_{n,h}^{\D}(\phi_t(K))$ is a decreasing function of $t \leq 0$. 
\end{theorem}

We obtain the following results concerning the asymptotic behavior of the hyperbolic area and the hyperbolic $n$-th diameter, as $t \to -\infty$. 
\begin{theorem}\label{harea}
Suppose $(\phi_t)$ is a non-elliptic semigroup of holomorphic functions in $\D$, $\Delta$ is a petal of $(\phi_t)$ and $K \subset \Delta$ compact. Then $$\lim_{t \to -\infty } \Ar^{\D}_h(\phi_t (K))= \Ar^{\Delta}_h K$$ and  $$ \lim_{t \to -\infty} d_{n,h}^{\D}(\phi_t (K))= d_{n,h}^{\Delta} K. $$
\end{theorem} 

Last but not least, we pursue to extend the results for condenser capacity of \cite{bkkp2} in the case of backward orbits. 
The ordered pair $(\D, \phi_t(K)) $ forms a condenser, as $\D$ is a subdomain of $\widehat{\mathbb{C}}$ and $\phi_t(K) $ is a compact subset of $\D$. In the case where $K$ is non-polar, then so is $\phi_t(K)$; \cite[Corollary 3.6.6]{ransford}. This allows us to measure the size of the condenser by means of its capacity; more information on condensers follows in Subsections \ref{diameter} and \ref{green}. 
We obtain the following result concerning the convergence of the capacity of $(\D, \phi_t(K)) $, as $t \to - \infty$. 

\begin{theorem}\label{capacity}
Suppose $(\phi_t)$ is a non-elliptic semigroup of holomorphic functions in $\D$, $\Delta$ is a petal of $(\phi_t)$ and $K \subset \Delta$ compact and non-polar. Then $$\lim_{t \to -\infty } \Capac(\D, \phi_t (K))= \Capac (\Delta, K). $$
\end{theorem}
 
The above results are restricted in the case where the backward orbit of $K$ is regular. The question that arises is how the characteristics of $\phi_t(K)$ change, provided the backward orbit of $K$ is non-regular.
In this case, we can talk about ``degenerate petals'', which are basically non-regular maximal invariant curves of $(\phi_t)$; we provide information on the variety of the forms petals can take in Subsection \ref{semigroups}.
We obtain the following outcome: 

\begin{theorem}\label{nonreg}
Let $(\phi_t)$ be a non-elliptic semigroup in $\mathbb{D}$. 
Suppose $\Delta$ is a degenerate petal of $(\phi_t)$ and $K$ is a compact non-polar subset of $\Delta$. Then
\begin{enumerate}[(i)]
    \item $\lim\limits_{t\to-\infty}\omega(\phi_t(z),\phi_t(K),\mathbb{D})=0$, for all $z\in \Delta \setminus K$,
    \item $\lim\limits_{t\to-\infty}\lambda_\mathbb{D}(\phi_t(z))=+\infty$, for all $z\in \Delta$,
    \item $\lim\limits_{t\to-\infty}d_\mathbb{D}(\phi_t(z),\phi_t(w))=+\infty$, for all $z,w\in \Delta$ with $z\neq w$,
    \item $\lim\limits_{t\to-\infty}g_\mathbb{D}(\phi_t(z),\phi_t(w))=0$, for all $z,w\in \Delta$ with $z\neq w$,
    \item $\lim\limits_{t\to-\infty}\Ar_{h}^\mathbb{D}(\phi_t(K))=0$,
    \item $\lim\limits_{t\to-\infty}d_{n,h}^{\D}(\phi_t(K))=+\infty$,
    \item $\lim\limits_{t\to-\infty}\Capac(\mathbb{D},\phi_t(K))=+\infty$.
\end{enumerate}
\end{theorem}

The article is structured as follows. In Section \ref{prepare}, we state some basic tools of one-parameter semigroups, potential theory and hyperbolic geometry, which will be used in proving the above Theorems. In Section \ref{harmea}, Theorems \ref{monotharm} and \ref{th2} are proved concerning the monotonicity and the asymptotic behavior of harmonic measure. 
Afterwards, the asymptotic behavior of hyperbolic metric, hyperbolic area and hyperbolic $n$-th diameter is examined in Section \ref{area&diam}, whereas in Section \ref{condenser}, similar results are obtained regarding the condenser capacity. Meanwhile in Section \ref{nonregular}, the case of a non-regular backward orbit of a compact set is investigated.  
 
\section{\bf Preparation for the proofs}\label{prepare}

\subsection{Koenigs Function-Petals-Backward Orbits}\label{semigroups}

For every non-elliptic one-parameter semigroup, there exists a Riemann mapping $h$ fixing the origin such that $\Omega: = h(\D) $ is a simply connected domain and also, \textit{convex in the positive direction}. This means that $\{ w+s: s\ge0 \} \subset \Omega$, for every $w \in \Omega$.
The function $h$ is unique up to a constant and it is called the \textit{Koenigs function} of the semigroup. A major property of the Koenigs function is that it linearizes the trajectories of the points in $\D$ under the semigroup; basically 
\begin{equation}\label{mapofh}
h(\phi_t(z) )= h(z)+t, 
\end{equation}
for all $z \in \D$ and $t \geq 0$.
The linearization property of the Koenigs function can also be generalized to the backward trajectories, supposing they exist. 

Before we move on to petals and the convergence of backward orbits, we need the following definition concerning the fixed points of a non-elliptic semigroup. 

\begin{definition}\cite[Chapters 12 \& 14]{Booksem}
\normalfont
Let $(\phi_t)$ be a non-elliptic semigroup. A boundary fixed point $\sigma \in \partial \D$ of $(\phi_t)$ is called \emph{regular}, if the angular derivative  $\phi_t^{\prime}(\sigma) < +\infty$. 

Suppose $\sigma \neq \tau$ is a boundary fixed point of $(\phi_t)$. 
If $\phi_t^{\prime}(\sigma)$ is finite, then $\sigma$ is a \emph{repelling fixed point} of $(\phi_t)$. Otherwise, $\sigma $ is a \emph{super-repelling fixed point} of $(\phi_t)$. 
A super-repelling fixed point is \emph{of first type} if it is the starting point of a maximal invariant curve of $(\phi_t)$.
\end{definition}

As stated in the Introduction in the case of regular backward orbits, for every point $z \in \Delta$, where $\Delta$ is a petal of $(\phi_t)$, the backward orbit exists and converges to a boundary regular fixed point. Hence it converges either to the Denjoy-Wolff point or to a repelling boundary fixed point of $(\phi_t)$.

In the case where $(\phi_t)$ is a group, the backward trajectory $\gamma_z $ is regular. If the group is hyperbolic and $\sigma \in \T \setminus \{ \tau \}$ is the repelling fixed point, then $\gamma_z (t)$ converges non-tangentially to $\sigma$, as $t \to - \infty$. On the other hand, if the group is parabolic, the backward orbit $\gamma_z (t)$ converges tangentially to the Denjoy-Wolff point $\tau$, as $t \to - \infty$. 

If $(\phi_t)$ is a not a group, for every backward orbit $\gamma_z$, there exists $\sigma \in \T$ (possible even $\sigma =\tau$) that is a fixed point of $(\phi_t)$ and $ \gamma_z(t) \xrightarrow{t \to -\infty} \sigma.$
In the case where $(\phi_t)$ is a hyperbolic semigroup (not a group) and $\sigma \in \T \setminus \{ \tau \}$ is a repelling fixed point, then $\gamma_z (t) \xrightarrow{t \to + \infty} \sigma$ non-tangentially.

\begin{figure}[h]
    \centering
    \includegraphics[scale=0.3]{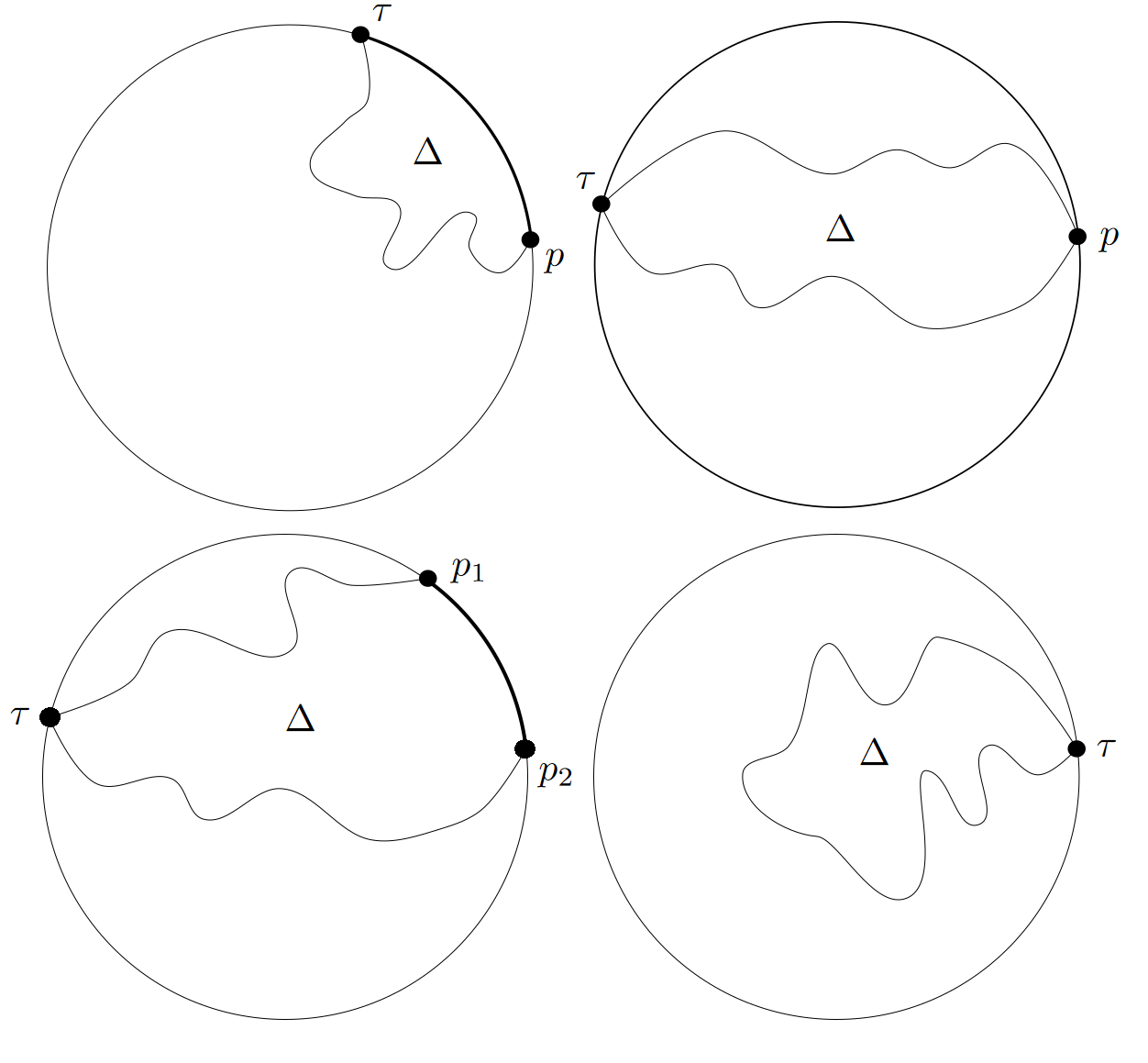}
    \caption{Hyperbolic \& Parabolic Petals}
    \label{typespetals}
\end{figure}

Furthermore, there exists the following characterization of petals of a non-elliptic semigroup, which is not a group. 
Let $\Delta$ be a petal of $(\phi_t)$.
The petal $\Delta$ is called \textit{hyperbolic} if $\partial \Delta$ contains a repelling fixed point of  $(\phi_t)$. If $\partial \Delta \setminus \{ \tau \}$ contains no boundary fixed points, then $\Delta$ is a \textit{parabolic petal}.
More specifically, only parabolic semigroups can have parabolic petals.
Furthermore, the boundary of $\Delta$ can contain at msugawaost two boundary fixed points; see \cite[Prop.13.4.10]{Booksem}. In Figure \ref{typespetals}, we observe all the possible cases on hyperbolic and parabolic petals. 

\begin{figure}[h]
    \centering
    \includegraphics[scale=0.4]{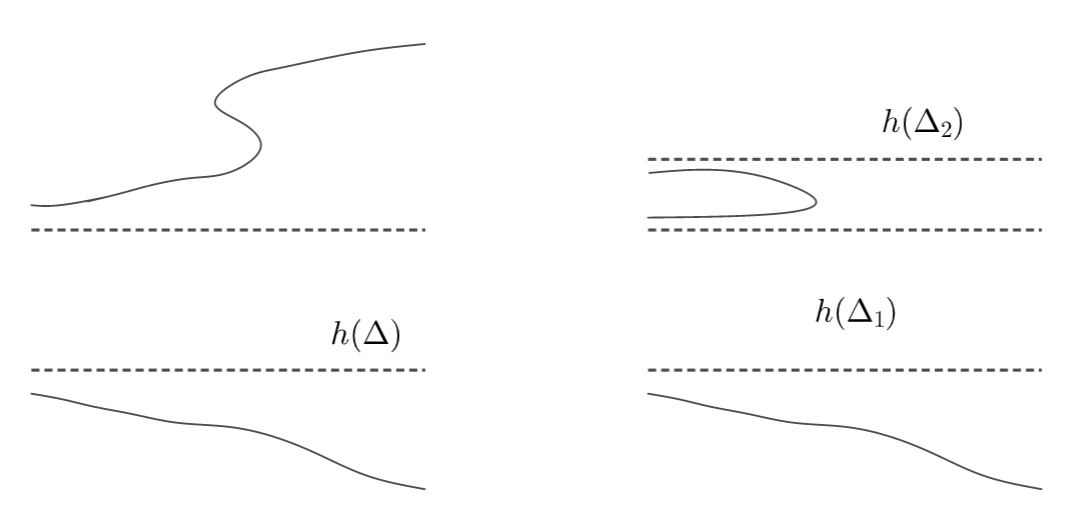}
    \caption{Images of Petals under $h$}
    \label{imagespetals}
\end{figure}

Let $h$ be the associated Koenigs function of $(\phi_t)$. In the case where $\Delta $ is a hyperbolic petal, then $h(\Delta)$ is a maximal horizontal strip contained in $\Omega$. 
When $\Delta$ is a parabolic petal, then $h(\Delta)$ is a maximal horizontal half-plane in $\Omega$; see Figure \ref{imagespetals}. 

As a result, petals have a direct connection with regular backward orbits. When it comes to non-regular backward orbits, it makes no sense to talk about petals in the way they were defined in the Introduction.  
A non-regular backward orbit can fall into one of the following three cases:
\begin{enumerate}[(i)]
    \item either it is part of the boundary of a hyperbolic petal, in which case it converges tangentially to a repelling fixed point of the semigroup
    \item either it is part of the boundary of a parabolic petal in which case it converges tangentially to the Denjoy-Wolff point of the semigroup
    \item or it converges to a super-repelling fixed point of the semigroup, in which case this super-repelling fixed point is of the first type and the convergence can be either tangential or non-tangential.
\end{enumerate}
Hence, in any of the above cases there exists a ``degenerate'' petal.
Let $\gamma:[0,+\infty)\to\mathbb{D}$ be a non-regular backward orbit for a non-elliptic semigroup $(\phi_t)$. Through the Koenigs function $h$, we move to the associated planar domain $\Omega$. Then, as with the regular backward orbits, the image $h(\gamma[0,+\infty))$ is a half-line that converges to $\infty$ through the negative direction. Remembering that the images of hyperbolic petals through $h$ are maximal horizontal strips in $\Omega$ and that the respective images of parabolic petals are maximal horizontal half-planes in $\Omega$, we can see that the image under $h$ of this degenerate petal is just the line that contains the set $h(\gamma[0,+\infty)).$ As a result, the image of such a petal under the Koenigs function is  $h(\Delta)=\{h(\gamma(0))+t:t\in\mathbb{R}\}$. 

For further information on the backward orbits in conjunction with the classification of one-parameter semigroups and the geometry of petals, the reader may refer to \cite[Chapter 13]{Booksem}.

\subsection{Hyperbolic and Quasi-hyperbolic metric}\label{hypmetric}

The \textit{hyperbolic metric} in $\D$ is $$ \lambda_{\D} (z) |dz|= \frac{|dz|}{1-|z|^2}, $$ where $\lambda_{\D}$ denotes its density.
If $U$ is a simply connected subdomain of the unit disk, then for $z \in U$, $\lambda_{\D}(z) \leq \lambda_{U} (z)$. 
The \textit{hyperbolic distance} between two points $a$, $b \in U$ is 
$$d_{U} (a,b)= \inf_{\gamma \subset U} \int_{\gamma} \lambda_{U} (z) |dz|,$$
where $\gamma$ is any rectifiable curve that lies in $U$ and joins $a$, $b$. The infimum is attained for the hyperbolic geodesic arc that joins $a,b$. 
For instance, the hyperbolic distance in the unit disk, for $z,w \in \D$ is equal to
$$ d_{\D} (z,w)= \arctanh \left| \frac{z-w}{1- \bar{z}w} \right| .$$
Hyperbolic distance is invariant under conformal mappings. Suppose that $f: \D \to \Omega$ is a conformal mapping, then $d_{\D} (z,w) = d_{\Omega} (f(z),f(w)),$ for every choice of $z,w \in \D$. 
In addition, if $K$ is a compact subset of $\D$, its hyperbolic area is given by the formula 
$$\Ar_h^{\D} (K)= \int_K \lambda_{\D}(z)^2 d A(z),$$
where $A$ is the Lebesgue area measure. 
Let's note that the hyperbolic area is also conformally invariant.
For two points $z_1, z_2$ in a simply connected domain $U$, the quasi-hyperbolic distance is defined as 
$$d^{\star}_U (z_1,z_2) = \min_C \int_C  \frac{|dw|}{\dist(w, \partial U)},$$ 
where the minimum is taken over all curves $C$ in $U$ joining $z_1$ and $z_2$; see e.g. \cite[p.92]{pommerenke4}. 

The following inequality connects hyperbolic and quasi-hyperbolic distances
\begin{equation}\label{relationquasidis}
\frac{1}{4} d^{\star}_U (z,w) \leq d_U (z,w) \leq d^{\star}_U (z,w),
\end{equation}
for every $z, w\in U$.
The reader may refer to \cite{bearmin, pommerenke4} for further properties of the hyperbolic metric. 

\subsection{Euclidean and Hyperbolic $n$-th Diameter - Capacities and Condensers} \label{diameter}
Let $K$ be a compact subset of $\mathbb{C}$. The \textit{Euclidean $n$-th diameter} of $K$ is 
\begin{equation}\label{eudiameter}
d_n (K) = \sup_{w_1,...,w_n \in K} \prod_{1 \leq \mu < \nu \leq n} |w_{\mu} - w_{\nu}|^{\frac{2}{n(n-1)}}
\end{equation}
and the supremum is attained, since $K$ is compact, for a $n$-tuple of points, which is called a \textit{Fekete $n$-tuple} for $K$; see \cite[Definition 5.5.1]{ransford}. We should point out that a Fekete $n$-tuple is not unique for the compact set $K$.
Its \textit{logarithmic capacity} $\Capac K$ is equal to the limit of $d_n(K)$, as $n \to + \infty$ . If two sets $A,B \subset \mathbb{C}$ differ on a set of zero logarithmic capacity (meaning $\Capac( A \setminus B)= \Capac (B \setminus A)=0$), we say that they are \textit{nearly everywhere} (n.e.) equal. Sets of zero logarithmic capacity are called \textit{polar sets} and they are negligible from the point of view of potential theory. 

Furthermore, for a compact set $K \subset \D$, its \textit{hyperbolic $n$-th diameter} is defined as 
$$ d_{n,h}^{\D} (K) = \sup_{w_1,...,w_n \in K} \prod_{1 \leq \mu < \nu \leq n} 
\left| \frac{w_{\mu} - w_{\nu}}{1 -  \overline{w}_{\mu}w_{\nu} } \right|^{\frac{2}{n(n-1)}},$$
where the supremum is attained for a $n$-tuple of points. 
The \textit{hyperbolic capacity} of $K$ is defined to be 
\begin{equation*}\label{hypcapacity}
\caph K = \lim_{ n \to + \infty} d_{n,h}^{\D} (K) 
\end{equation*}
and it is a conformally invariant quantity, due to the conformal invariance of the hyperbolic distance.
 
Another potential theoretic and conformally invariant quantity is the capacity of a condenser. 
A \textit{condenser} is an ordered pair $(D, K)$, where $D$ is a proper domain of $\widehat{\mathbb{C}}$ and $K $ is a compact subset of $D$. 
Suppose both $\partial D$ and $K$ are non-polar. In the special case where $D$ is a simply connected domain, the \textit{capacity} of $(D, K)$ can be defined as 
$$ \Capac (D, K) :=  2 \pi \caph_D K,$$
where $\caph_D K$ denotes the hyperbolic capacity of $K$ in $D$. 

As stated above, the condenser capacity is conformally invariant. 
Suppose $D$ and $G$ are simply connected domains in $\mathbb{C}$. If $f:D \to G$ is holomorphic and $K\subset D$, then 
\begin{equation}\label{confinvcapacity}
\Capac(D,K) \geq  \Capac (G,f(K)). 
\end{equation}
Equality holds if and only if $f$ is conformal. For more information on condensers, the reader may refer to \cite{dubininbook}.

\subsection{Green Potential - Green capacity}\label{green}
Let $D$ be a domain of the extended complex plane $\widehat{\mathbb{C}}$. 
The \textit{Green function} of $D$ with \textit{pole} at $w \in D$ is denoted by $g_D(\cdot, w)$ and satisfies the following conditions:
\begin{enumerate}[(i)]
\item $g_D(\cdot,w) $ is positive harmonic on $D \setminus \{w \}$ and bounded outside every neighborhood of $w$
\item $g_D(\cdot,w)+ \log| \cdot -w|$ is harmonic on $D$ 
\item $g_D(w,w)=\infty$, and for $z \to w$
$$ g_D(z,w) =   \begin{cases}
 \log |z| + \mathcal{O} (1), \, &w=\infty \\
	-\log|z-w| + \mathcal{O}(1), \, &w \neq \infty
  \end{cases} $$
\item $g_D(\cdot,w) = 0$ on the boundary $ \partial D$. 
\end{enumerate}
For every simply connected domain $D$, there exists the following relation between its Green function and the hyperbolic distance:  
\begin{equation}\label{greenformula}
g_{D}(z,w) = \log \tanh d_{D}(z,w),
\end{equation}
for $z, w \in \D$; see \cite[p.109]{ransford}. 
If the boundary of a domain $D$ is non-polar, then the Green function $g_{D}$ exists and it is unique. In this case, the domain $D$ is called \textit{Greenian}. The Green function is symmetric,
for every $z,w \in D$. Moreover, the Green function is conformally invariant; e.g. \cite[Theorem 4.4.4]{ransford}.

Let $D$ be a Greenian domain of $\widehat{\mathbb{C}}$ with Green function $g_D (x,y)$, for $x,y \in D$. For a measure $\mu$ with compact support in $D$, its \textit{Green energy} is defined as the integral 
$$I_D[\mu] : = \iint g_D(x,y) d \mu(x) d\mu(y) .$$

Suppose $K$ is a compact subset of $D$. 
The \textit{Green energy} of $K$ with respect to $D$ is defined as 
$$V(K,D) = \inf_{\mu} I_D[\mu],$$
where the infimum is taken over all the Borel measures $\mu$ with compact support $K$ such that $\mu(K)=1$. If $V(K,D)< + \infty$, the infimum is attained for a Borel measure $\mu$, which is called \textit{Green equilibrium measure} of $K$.
The Green energy of a compact set is directly associated with condenser capacity. 
\begin{remark}
If $K$ and $\partial D$ have positive logarithmic capacity, the capacity of the condenser $(D,K)$ is proportional to the Green energy of the compact set $K$ and it is true that
\begin{equation}\label{greencapacity2}
\Capac (D, K)  = \frac{2 \pi }{V(K,D)}.
\end{equation}
\end{remark}

More information on Green energy and the aspects of Green function can be found in \cite{gardiner, helms}.

\subsection{Harmonic Measure}\label{harmonicmeasure}

Let $D$ be a proper subdomain of $\widehat{\mathbb{C}}$ with non-polar boundary and $\mathcal{B}(\partial D)$ be the $\sigma$-algebra of the Borel sets of $\partial D$.   
Suppose $E \in \mathcal{B}(\partial D)$. The \textit{harmonic measure} of $E$ at a point $z \in D$ is the solution of the generalized Dirichlet problem in $D$ with boundary values $1$ on $E$ and $0$ on $\partial D \setminus E$. 
For a fixed $E \in \mathcal{B}(\partial D)$, $\omega ( \cdot , E, D)$ is a harmonic and bounded function on $D$. 
Moreover, for a fixed point $z \in D$, the map 
$$ \omega ( z, \cdot, D) : \mathcal{B}(\partial D)  \to [0,1] \quad \text{with} \quad 
E \mapsto \omega ( z, E , D)$$
is a Borel probability measure on $\partial D$. 
In addition, if $\zeta$ is a regular point of $\partial D$ which lies outside the relative boundary of $E$ in $\partial D$, then $ \lim\limits_{ z \to \zeta} \omega ( z, E, D)= \mathcal{X}_E (\zeta)$, where $\mathcal{X}_E (\cdot) $ denotes the characteristic function of $E$.

A major property of the harmonic measure is its conformal invariance; see \cite[\S\ 4.3]{ransford}.
For the sake of simplicity, if $E$ is a compact subset of $D$ with positive logarithmic capacity, from now on we will use the notation $\omega(z, E, D) := \omega(z, \partial E, D \setminus E).$ Let us state the following property for the harmonic measure.

\begin{proposition}[Strong Markov Property for Harmonic Measure]\cite[p.88]{portstone}\label{markovharmonic}
Suppose $\Omega$ is a Greenian domain on $\mathbb{C}$ and $S$ a domain such that $S \subset \Omega$. Let $E\subset \partial \Omega$. Set $A:= \partial S \cap \Omega$. Then for $z \in \Omega $, 
$$\omega(z, E, \Omega) = \omega (z, E, S)+ \int_A \omega( s, E, \Omega)\cdot \omega (z, d s, S). $$
\end{proposition} 

Respectively, we state the following relation between harmonic measure and the Green function. 

\begin{proposition}[Strong Markov Property for Green function]\cite[p.111]{portstone}\label{markovgreen}
Suppose $\Omega$ is a Greenian domain on $\mathbb{C}$ and $S$ a domain such that $S \subset \Omega$. Set $A:= \partial S \cap \Omega$. Then for $z,w \in S$, 
$$g_{\Omega}(z,w) =g_S(z,w) + \int_A g_{\Omega}(\alpha , z)\cdot \omega (w, d \alpha, S). $$
\end{proposition} 

Another property of great significance for harmonic measure is its probabilistic interpretation. Suppose $D$ is a domain on $\mathbb{C}$ and $E$ a Borel subset of $\partial D$. Let $B_t$, $t>0$, be a Brownian motion on the complex plane starting from a point $z \in D$. Let $t_0= \inf \{ t>0 : B_t \notin D \}$ be the first exit time of $B_t$ from $ D$. The harmonic measure $\omega(z, E, D)$ is the probability of $B_{t_0} \in E$. 
Information and detailed theory on harmonic measure can be found in \cite{margarnett} and \cite[Chapter 4]{ransford}. 

\section{ \bf Harmonic Measure - Proof of Theorems \ref{monotharm} and \ref{th2}}\label{harmea}

In the course of the following paragraphs, we suppose that $(\phi_t)$ is a non-elliptic semigroup with associated Koenigs function $h$ and $\Delta$ is a petal of $(\phi_t)$. Suppose $\Omega$ is the associated Koenigs domain of the semigroup. Take $K$ to be a compact non-polar subset of $\Delta$. 

\textbf{Proof of Theorem \ref{monotharm}.} Fix a point $z \in \Delta \setminus K$. Set $s \in \mathbb{R} $ such that $t<s \leq 0$. According to the subordination principle of harmonic measure (see \cite[Theorem 4.3.8]{ransford}), 
\begin{eqnarray*}
\omega (\phi_s (z), \phi_s (K), \D)& =&  \omega (h(z)+s, h(K)+s, \Omega) \\
&=& \omega (h(z), h(K), \Omega -s) \\
&\geq & \omega (h(z), h(K), \Omega -t) \\ 
&=& \omega (\phi_t(z), \phi_t (K), \D), 
\end{eqnarray*}
as $\Omega -t \subset \Omega -s$. Therefore Theorem \ref{monotharm} is proved. 
\qed

\textbf{Proof of Theorem \ref{th2}.} Fix a point $z \in \Delta \setminus K$. Using the monotonicity property of the harmonic measure and considering consecutively the fact that $\phi_t$ is an automorphism of $\Delta$, for all $t\leq0$, we see that
   \begin{eqnarray*}
   \omega(\phi_t(z),\phi_t(K),\mathbb{D})&\ge&\omega(\phi_t(z),\phi_t(K),\Delta)\\
   &=&\omega(z,K,\Delta),
   \end{eqnarray*}
   for all $t\leq0$. Therefore,
   \begin{equation}
       \lim\limits_{t\to-\infty}\omega(\phi_t(z),\phi_t(K),\mathbb{D})\ge\omega(z,K,\Delta),
   \end{equation}
  since the limit exists due to the monotonicity established in Theorem \ref{monotharm}. 
Now, we need a similar reverse inequality. Due to conformal invariance and the properties of the Koenigs function, for $t \leq 0$ and $z\in\Delta\setminus K$, $$\omega(\phi_t(z),\phi_t(K),\mathbb{D})=\omega(h(\phi_t(z)),h(\phi_t(K)),\Omega)=\omega(h(z)+t,h(K)+t,\Omega). $$
Since $K$ is compact, there exists a $t_0\leq0$ such that $h(z)+t\notin h(K)$, for all $t\leq t_0$. 
We will deal first with the case when the petal $\Delta$ is hyperbolic. Then $h(\Delta)$ is a maximal strip contained in $\Omega$. 

\begin{figure}[h]
       \centering
       \includegraphics[scale=0.52]{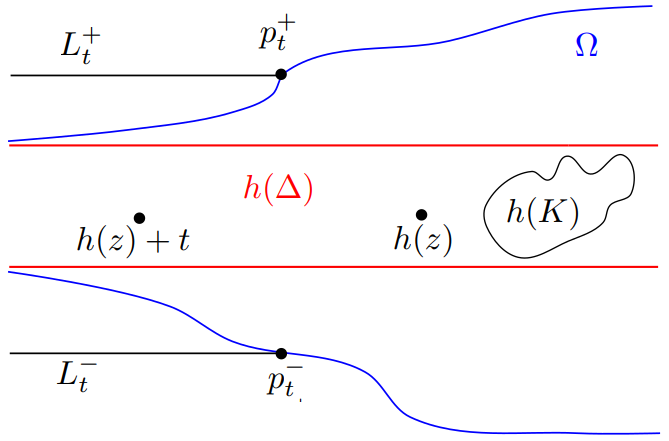}
       \caption{The construction of $\Omega_t$}
       \label{constructon1.2}
   \end{figure}

In combination with the fact that $\Omega$ is convex in the positive direction, there exists $t_1\le0$ such that the sets 
$$\partial\Omega\cap\{\zeta:\IM\zeta>\IM h(z), \RE\zeta=\RE h(z)+t\} \quad  \text {and} \quad \partial\Omega\cap\{\zeta:\IM\zeta<\IM h(z), \RE\zeta=\RE h(z)+t\}$$ are both non-empty, for all $t\le t_1$. 
For $t<2\min\{t_0,t_1\}$, we denote by $p_t^+$ the point of $\partial\Omega$ such that $$\left|p_t^+-\left(h(z)+\frac{t}{2}\right)\right|=\min\left\{\left|\zeta-\left(h(z)+\frac{t}{2}\right)\right|:\zeta\in\partial\Omega, \IM\zeta>\IM h(z), \RE\zeta=\RE h(z)+\frac{t}{2}\right\}.$$ 
Similarly, we denote by $p_t^-$ the point of $\partial\Omega$ such that $$\left|p_t^- -\left(h(z)+\frac{t}{2}\right)\right|=\min\left\{\left|\zeta-\left(h(z)+\frac{t}{2}\right)\right|:\zeta\in\partial\Omega, \IM\zeta<\IM h(z), \RE\zeta=\RE h(z)+\frac{t}{2}\right\}.$$ 
Then we define the sets $$L_t^+ : =\left\{\zeta:\RE\zeta\le\RE p_t^+, \IM\zeta=\IM p_t^+\right\} \quad \text{and} \quad  L_t^- : =\left\{\zeta:\RE\zeta\le\RE p_t^-, \IM\zeta=\IM p_t^-\right\}$$ and we construct the domain $\Omega_t :=\mathbb{C}\setminus(L_t^+\cup L_t^-)$ (see Figure \ref{constructon1.2}).
   
   Since the families of points $(p_t^+)_t,(p_t^-)_t\subset\partial\Omega$, by convexity it is clear that $\Omega\subset\Omega_t$, for all sufficiently small $t$. According to the domain monotonicity of the harmonic measure, we obtain 
\begin{eqnarray*}
   \omega(h(z)+t,h(K)+t,\Omega)&\le&\omega(h(z)+t,h(K)+t,\Omega_t)\\&=&\omega(h(z),h(K),\Omega_t-t).
   \end{eqnarray*}
   It is easy to see that $h(\Delta)\subset\Omega_t-t$, for all sufficiently small $t$. Moreover, $$\partial h(K)\subset\partial(h(\Delta)\setminus h(K))\cap\partial((\Omega_t-t)\setminus h(K))$$ and $$\partial(h(\Delta)\setminus h(K))\cap((\Omega_t-t)\setminus h(K))=\partial h(\Delta).$$ 
   Since $z\in\Delta\setminus K$, we have $h(z)\in h(\Delta)\setminus h(K)$ and applying Proposition \ref{markovharmonic}, we have
   $$\omega(h(z),h(K),\Omega_t-t)=\omega(h(z),h(K),h(\Delta))+\int\limits_{\partial h(\Delta)}\omega(\zeta,h(K),\Omega_t-t)\cdot\omega(h(z),d\zeta,h(\Delta)\setminus h(K)).$$
   By the construction above, we can see that 
   \begin{equation}\label{realpt}
     \lim\limits_{t\to-\infty}(\RE p_t^+-t)=\lim\limits_{t\to-\infty}(\RE p_t^--t)=+\infty,  
   \end{equation}
   while for $h(\Delta)=\{ w :a<\IM w<b\}$, $$\lim\limits_{t\to-\infty}\IM p_t^+=b \quad \text{and} \quad \lim\limits_{t\to-\infty}\IM p_t^-=a, $$ 
   since the strip $h(\Delta)$ is maximal in $\Omega$. Hence the boundary of $\Omega_t -t$, as $t$ decreases, approaches $\partial h(\Delta)$. 
   We will first prove that for any $\zeta\in\partial h(\Delta)$, we have $ \omega(\zeta,h(K),\Omega_t-t)\xrightarrow{t\to-\infty} 0.$
   
   The set $h(K)$ is compact and so, there exists a horizontal strip $S$ such that $h(K) \subset S \subsetneq h(\Delta)$. We denote by $\partial S^{+}$ and $\partial S^{-}$ the upper and lower boundary components of $S$, respectively.  
   Set $A_t^{+}$ and $A_t^{-}$ to be the horizontal lines that contain $L_t^{+}$ and $L_t^{-}$, respectively. 
   Define $S_1^{t}$ to be the horizontal strip bounded by $A_t^{+}$ and $\partial S^{+}$ and $S_2^{t}$ the one bounded by $\partial S^{-}$ and $A_t^{-}$. 
  Without loss of generality, we take a point $\zeta \in \partial h(\Delta) \cap S_1^{t}$. If $\zeta \in S_2^{t}$, the proof follows in the same manner. We write $S_1^t = \{ w \in \mathbb{C} : a_1< \IM w < \IM p_t^+ \}$. According to \cite[p. 100]{ransford}, the harmonic measure 
  $$\omega( \zeta, A_t^{+}, S_1^t)= \frac{\IM \zeta -a_1}{\IM p_t^{+} - a_1} = \frac{b -a_1}{\IM p_t^{+} - a_1}.$$ However, $\IM p_t^{+} \to b$, as $t \to - \infty$, from the construction of $\Omega_t$'s. As a result, 
  $ \omega( \zeta, A_t^{+}, S_1^t) \xrightarrow{t \to -\infty} 1 $. 
  Furthermore, from \eqref{realpt}, we obtain that the half-line $L_t^{+}-t$ expands towards $\infty$, in the positive direction. Hence, when $ t \to -\infty$, $L_t^{+}-t$ and $A_t^{+}$ coincide and so, $$\lim_{t \to - \infty} \left[ \omega( \zeta, A_t^{+}, S_1^t) - \omega(\zeta, L_t^{+}-t, S_1^t) \right]=0$$ that implies $\omega(\zeta, L_t^{+}-t, S_1^t) \xrightarrow{ t \to -\infty} 1$. However $S_1^t \subset (\Omega_t -t) \setminus h(K)$ and according to the subordination principle of harmonic measure, 
  $$\omega(\zeta, L_t^{+}-t, S_1^t) \leq \omega (\zeta , L_t^{+}-t, (\Omega_t -t) \setminus h(K))\Rightarrow \lim_{t \to -\infty} \omega (\zeta , L_t^{+}-t, (\Omega_t -t) \setminus h(K)) \geq 1. $$
  From the properties of harmonic measure, $\omega (\zeta , L_t^{+}-t, (\Omega_t -t) \setminus h(K)) \xrightarrow{t \to -\infty} 1$ and we deduce that $$ \lim_{t \to -\infty} \omega (\zeta, h(K), \Omega_t-t) =0.$$ The choice of $\zeta$ was arbitrary and hence, the convergence is true for all $\zeta\in\partial h(\Delta)$. 
  Next, we examine the uniform convergence. Consider the quantity $$\sup\limits_{\zeta\in\partial h(\Delta)}\omega(\zeta,h(K),\Omega_t-t).$$ Since $h(K)$ is compact, we have $$\lim\limits_{ \Omega_t -t \, \ni \zeta \to\infty}\omega(\zeta,h(K),\Omega_t-t)=0.$$
    Therefore, the above supremum is attained on some point of $\partial h(\Delta)$. We denote by $\zeta_t$ such a point and we distinguish two cases: either $\RE \zeta_t>\RE p_t^+ -t$, for all $t<2\min\{t_0,t_1\}$ or $\RE \zeta_t\le\RE p_t^+ -t$, for all $t<2\min\{t_0,t_1\}$. Of course, these two cases might alternate as $t\to-\infty$, but, in such a scenario, we can just consider two subfamilies of $(\zeta_t)$ and proceed in the same manner. 
   In the first case, it is clear that $\zeta_t \xrightarrow{t \to -\infty} \infty$, something which directly implies that $\omega(\zeta_t,h(K),\Omega_t-t) \xrightarrow{t\to-\infty}0.$
   In the second case, as $t\to-\infty$, one of the half-lines $L_t^+-t,L_t^--t$, and thus the boundary of $\Omega_t -t$, is getting arbitrarily close to $\zeta_t$. 
   Using the above construction of horizontal strips, we are led again to  $\omega(\zeta_t,h(K),\Omega_t-t) \xrightarrow{t\to-\infty}0.$ 
Therefore, we get $$\sup\limits_{\zeta\in\partial h(\Delta)}\omega(\zeta,h(K),\Omega_t-t) \xrightarrow{t\to-\infty} 0$$
 and the harmonic measure converges uniformly to $0$ on $\partial h(\Delta)$. 
   Let $\epsilon>0$. Then, there exists $t_2<2\min\{t_0,t_1\}$ such that $\omega(\zeta,h(K),\Omega_t-t)<\epsilon,$ for all $\zeta\in\partial h(\Delta)$ and all $t\le t_2$. Returning to the Strong Markov Property, for $t\le t_2$, we get
   \begin{eqnarray*}
   \omega(h(z),h(K),\Omega_t-t)&<&\omega(h(z),h(K),h(\Delta))+\int\limits_{\partial h(\Delta)}\epsilon\cdot\omega(h(z),d\zeta,h(\Delta)\setminus h(K))\\
   &=&\omega(h(z),h(K),h(\Delta))+\epsilon\cdot\omega(h(z),\partial h(\Delta),h(\Delta)\setminus h(K))\\
   &\le&\omega(h(z),h(K),h(\Delta))+\epsilon,
   \end{eqnarray*}
   since the harmonic measure attains values in $[0,1]$. Therefore, 
   $$\limsup\limits_{t\to-\infty}\omega(h(z),h(K),\Omega_t-t)\le\omega(h(z),h(K),h(\Delta))=\omega(z,K,\Delta).$$
   With the use of conformal invariance of the harmonic measure, we obtain 
   \begin{eqnarray*}
   \lim\limits_{t\to-\infty}\omega(\phi_t(z),\phi_t(K),\mathbb{D})&=&\lim\limits_{t\to-\infty}\omega(h(z)+t,h(K)+t,\Omega)\\
   &\le&\limsup\limits_{t\to-\infty}\omega(h(z),h(K),\Omega_t-t)\\
   &\le&\omega(z,K,\Delta).
   \end{eqnarray*}
   All in all, 
   $$\lim\limits_{t\to-\infty}\omega(\phi_t(z),\phi_t(K),\mathbb{D}) =\omega(z,K,\Delta)$$
   and we have the desired result.
   
    If the petal $\Delta$ is parabolic, then $h(\Delta)$ is a horizontal half-plane. In this case, only one of the families $(p_t^+),(p_t^-)$ exists. Without loss of generality, we can assume that it is $(p_t^+)$. Then, we construct the half-lines $(L_t^+)$ in the same manner as before and consider the simply connected domains $\Omega_t=\mathbb{C}\setminus L_t^+$, for all suitable $t\le0$. 
    The set $h(K)$ is compact and so, there exists a horizontal half-plane $H$ such that $h(K) \subset H \subsetneq h(\Delta)$. The horizontal line $\partial H $ along with $A_t$, which is also a horizontal line that contains $L_t^+$, are the boundary components of a horizontal strip that contains $\partial h(\Delta)$. 
    Following the same procedure as in the hyperbolic case, we obtain that the harmonic measure in $\Omega_t-t \setminus h(K)$ with respect to $\partial h(K)$ converges uniformly to $0$ on $\partial  h (\Delta)$. Then, we continue with the proof exactly as above and deduce the desired result.
   \qed

\section{\bf Hyperbolic Geometric Quantities - Proof of Theorems \ref{limdist}, \ref{monothypdiam} and \ref{harea}}\label{area&diam}

\textbf{Proof of Theorem \ref{limdist}}. 
     Suppose $h$ is the associated Koenigs function of the semigroup $(\phi_t)$. Let $z\in \Delta$ with $w:=h(z)$. We recall that $(\phi_{t|\Delta})$ is a group of automorphisms of $\Delta$. Hence
$$\lambda_{\D}(\phi_t(z))  \leq \lambda_{\Delta} (\phi_t(z)) =  \lambda_{\Delta} (z), $$
where the inequality holds due to the domain monotonicity property of the hyperbolic density. 
Due to conformal invariance, the above inequality can be written as 
$$\lambda_{\Omega}(w+t) = \lambda_{\Omega -t} (w) \leq \lambda_{h(\Delta)}(w). $$
    
For $t \leq 0$, the family of domains $(\Omega -t )_t$ decreases, as $t$ decreases. 
According to \cite[Theorem 1]{mindaquotients}, since $h(\Delta) \subset \Omega -t$, we obtain
\begin{equation}\label{doubleeq}
    1 \leq \frac{\lambda_{h(\Delta)}(w)}{\lambda_{\Omega - t}(w)}\leq 1+ \frac{2}{e^{R_t(w)}-1}, 
\end{equation}
where 
\begin{equation}\label{distanceRt}
    R_t(w)= d_{\Omega - t}(w,(\Omega - t) \setminus h(\Delta)) = d_{\Omega -t}(w, \partial h(\Delta)) = \inf_{\zeta \in \partial h(\Delta)} d_{\Omega - t} (w, \zeta).
\end{equation}

Obviously, $$\lim_{\partial h(\Delta) \ni \zeta \to \infty }d_{\Omega-t}(w, \zeta)=+ \infty, $$ and hence the infimum in \eqref{distanceRt} is attained for a specific point on $\partial h(\Delta)$. We define the family of points $(\zeta_t)_t$ such that for every $t \leq 0$, 
\begin{equation}\label{infdistance}
    d_{\Omega-t} (w, \zeta_t) = \inf_{\zeta \in \partial h (\Delta)} d_{\Omega -t}(w, \zeta)= R_t(w).
\end{equation}
However, using the quasi-hyperbolic metric, we have that for every $\zeta \in \partial h(\Delta)$, 
$$\lambda_{\Omega-t}^{\star}(\zeta)= \frac{1}{\dist(\zeta, \partial \Omega-t)} \xrightarrow{t \to -\infty} +\infty, $$
due to the convergence of the family $(\Omega -t)_t$. Hence, for every $\zeta \in \partial h(\Delta)$, $d^{\star}_{\Omega -t} ( w, \zeta) \xrightarrow{t \to -\infty} + \infty$. 
Returning back to \eqref{infdistance} and using the connection between hyperbolic and quasi-hyperbolic distance \eqref{relationquasidis}, we obtain
\begin{eqnarray}\label{Rtlimit}
\lim_{t \to -\infty} R_t(w) &=& \lim_{t \to -\infty} d_{\Omega - t}(w,\zeta_t)  \nonumber \\
&\geq & \frac{1}{4}  \lim_{t \to -\infty} d^{\star}_{\Omega -t} ( w, \zeta_t) \nonumber\\
&=& + \infty.
\end{eqnarray}
Hence $$1\leq \lim_{t \to -\infty} \frac{\lambda_{h(\Delta)}(w)}{\lambda_{\Omega-t}(w)} \leq 1+ \lim_{t \to -\infty} \frac{2}{e^{R_t(w)}-1}=1,$$
which leads to 
\begin{equation}\label{convofdensity}
    \lim_{t \to -\infty} \lambda_{\Omega-t} (w) = \lim_{t \to -\infty} \lambda_{\Omega}(w+t) = \lambda_{h(\Delta)}(w).
\end{equation}
Now we examine the uniform convergence. We can rewrite \eqref{doubleeq} in the following way:
$$0 \leq \lambda_{h(\Delta)}(w) - \lambda_{\Omega - t}(w) \leq \frac{2\lambda_{\Omega - t}(w)}{e^{R_t(w)}-1}.$$
The hyperbolic density $\lambda_{\Omega -t}$ is bounded on compacta in $h(\Delta)$. Suppose $a$ lies on a compact subset of $\Delta$. Then, according to \eqref{Rtlimit}, for every $\epsilon>0$, there exists $t_1 \leq 0$ and $M >0$ such that $R_t(h(a)) > M$, for every $t \geq t_1$. Hence, there exists a constant $c>0$ such that 
$$0 \leq \lambda_{h(\Delta)}(w) - \lambda_{\Omega - t}(w) \leq \frac{2c}{e^M -1}. $$
However $M$ is arbitrarily large, so we let $M \to +\infty$ and it follows that $\lambda_{\Omega - t}(w)$ converges locally uniformly to $\lambda_{h(\Delta)}(w) $ in $h(\Delta)$. Due to conformal invariance, $\lambda_{\D}(\phi_t(\cdot))$ converges locally uniformly to $\lambda_{\Delta}(\cdot)$ in $\Delta. $ 

Next, we move on to the convergence of the hyperbolic distance. Suppose $z, w \in \Delta$. 
Due to the monotonicity property of the hyperbolic distance and the fact that every $\phi_t$, for $t\le0$, is an automorphism of $\Delta$, we have
     \begin{equation*}
     d_{\D}(\phi_t(z),\phi_t(w))\leq d_{\Delta}(\phi_t(z),\phi_t(w)) =d_\Delta(z,w),
      \end{equation*}
    for all $t\le0$. Therefore,
     \begin{equation}\label{limsup}
         \limsup\limits_{t\to-\infty}d_{\D}(\phi_t(z),\phi_t(w))\le d_\Delta(z,w).
     \end{equation}
    Since $h(\Delta)\subset\Omega-t$ for all $t\le0$, it is also true that $h(z),h(w)\in\Omega-t$, for all $t\le0$. So, let $\delta_t:[0,1]\to\Omega-t$ be the hyperbolic geodesic arc of $\Omega-t$ such that $\delta_t(0)=h(z)$ and $\delta_t(1)=h(w)$. Suppose that there exists a sequence of real numbers $\{t_n\}\subset(-\infty,0]$ with $\lim\limits_{n\to+\infty}t_n=-\infty$, such that $\delta_{t_n}[0,1]\cap(\Omega-t)\setminus h(\Delta)\neq\emptyset$. As a result, for any $n\in\mathbb{N}$, there exists at least one $x_n\in(0,1)$ such that $\delta_{t_n}(x_n)\in\Omega-t\setminus h(\Delta)$. Due to the conformal invariance of the hyperbolic distance and the fact that $\delta_{t_n}$ is a geodesic arc, we have
    \begin{eqnarray*}
    d_\Delta(z,w)&=&d_{h(\Delta)}(h(z),h(w))\\
    &\ge&d_{\Omega-t_n}(h(z),h(w))\\
    &\ge&d_{\Omega-t_n}(h(z),\delta_{t_n}(x_n))\\
    &\ge&\frac{1}{4}\log \left(1+\frac{|\delta_{t_n}(x_n)-h(z)|}{\min\{\delta_{\Omega-t_n}(h(z)),\delta_{\Omega-t_n}(\delta_{t_n}(x_n))\}}\right),
    \end{eqnarray*}
    where $\delta$ denotes the Euclidean distance from the boundary and the last inequality follows from \cite[Theorem 5.3.1]{Booksem}.
        We distinguish the following two cases:
    \begin{enumerate}[(i)]
        \item either $\inf\RE(\delta_{t_n}(x_n))\in\mathbb{R}$ and $\sup\RE(\delta_{t_n}(x_n))=+\infty$,
        \item or both $\inf\RE(\delta_{t_n}(x_n)),\sup\RE(\delta_{t_n}(x_n))\in\mathbb{R}$.
    \end{enumerate}
    Any other case where $\inf\RE(\delta_{t_n}(x_n))=-\infty$ can be treated in a similar manner as (i).
 In case (i), it is clear that $\delta_{t_n}(x_n)\xrightarrow{n\to+\infty}\infty$, while $\delta_{\Omega-t_n}(h(z))$ remains bounded from above. Therefore, by the inequalities above, we are led to $\delta_\Delta(z,w)=+\infty$. Contradiction! In case (ii), we can see that $\delta_{\Omega-t_n}(\delta_{t_n}(x_n))\xrightarrow{n\to+\infty}0$ (because $t_n\xrightarrow{n\to+\infty}-\infty$), while $|\delta_{t_n}(x_n)-h(z)|$ remains bounded from below. Therefore, once again, we are led to $\delta_\Delta(z,w)=+\infty$. Contradiction!
    
     As a result, there can be no such sequence $\{t_n\}$, which means that there exists a $T\le0$ such that $\delta_t([0,1])\subset h(\Delta)$, for all $t\le T$.
     Bearing in mind the uniform convergence of the hyperbolic density on compacta, for every $\epsilon>0$, there exists $t_0 \leq T$ such that, for all $t \leq t_0$,
     \begin{eqnarray*}
     d_{\Delta}(z,w) &=& d_{h(\Delta)} (h(z), h(w)) \\
     &\leq & \int_{\delta_t} \lambda_{h(\Delta) }(\zeta) |d\zeta| \\
     &<& \int_{\delta_t} \left[\lambda_{\Omega-t }(\zeta) +\epsilon\right] |d\zeta|\\
     &= & d_{\Omega -t} (h(z), h(w)) + \epsilon \cdot \len (\delta_t),
     \end{eqnarray*}
     where $\len(\delta_t)$ is the euclidean length of $\delta_t$. The curve $\delta_t$ is the geodesic joining $h(z)$ and $h(w)$ in $\Omega -t$ and so, its pre-image under the Koenigs function is the geodesic joining $\phi_t(z)$ and $\phi_t (w)$ in $\D$. We denote by $C$ the line segment $[h(z), h(w)]$ in $\Omega -t$. According to the Gehring-Hayman Theorem \cite[\S\ 4.1- 4.6]{pommerenke4}, there exists an absolute constant $K$ such that 
     \begin{eqnarray*}
      \len (\delta_t) & \leq & K \len (C) \\
      &=& K |h(z)-h(w)| < + \infty. 
      \end{eqnarray*}
    Therefore, we obtain that $$\liminf_{t \to -\infty} d_{\Omega -t} (h(z), h(w)) \geq d_{h(\Delta)} (h(z), h(w)), $$
    which combined with \eqref{limsup} and the conformal invariance of the hyperbolic metric, gives us the desired result. \qed
  
\textbf{Proof of Theorem \ref{monothypdiam}}. 
Fix $t,s \leq  0$. The hyperbolic $n$-th diameter $$d_{n,h}^{\D}(\phi_t(K))^{\frac{n(n-1)}{2}} = \max_{ z_1,..., z_n \in \phi_t(K)} \prod_{1 \leq \mu < \nu \leq n}  \tanh  d_{\D}( z_{\mu}, z_{\nu}) =\prod_{1 \leq \mu < \nu \leq n}  \tanh  d _{\D}( y_{\mu},y_{\nu}),$$
where $y_1, ..., y_n$ is a Fekete $n$-tuple of $\phi_t(K)$. The hyperbolic $n$-th diameter of $\phi_{t+s} (K)$ is equal to 
$$d_{n,h}^{\D}(\phi_{t+s}(K))^{\frac{n(n-1)}{2}}=\max_{ z_1,..., z_n \in \phi_t(K)} \prod_{1 \leq \mu < \nu \leq n}  \tanh  d_{\D}( \phi_s(z_{\mu}), \phi_s( z_{\nu})). $$
Using Schwarz-Pick lemma, we get
\begin{eqnarray*}
d_{n,h}^{\D}( \phi_{t}(K))^{\frac{n(n-1)}{2}} & = &\prod_{1 \leq \mu < \nu \leq n}  \tanh  d _{\D}( y_{\mu},y_{\nu})\\
& \leq &\prod_{1 \leq \mu < \nu \leq n} \tanh d_{\D}( \phi_s(y_{\mu}), \phi_s( y_{\nu}))  \\
&\leq & d_{n,h}^{\D} (\phi_{t+s} (K))^{\frac{n(n-1)}{2}},
\end{eqnarray*} 
for every choice of $t$. 
Therefore, $d_{n,h}^{\D} (\phi_t (K))$ is a decreasing function of $t \leq 0$. \qed

\textbf{Proof of Theorem \ref{harea}}. 
The hyperbolic area of $\phi_t(K)$ in $\D$ is 
\begin{eqnarray*}
\Ar_h^{\D} (\phi_t(K)) &=&\int_{\phi_t(K)}  \lambda_{\D} (z)^2 d A(z) \\
&=&\int_{h(K)} \lambda_{\Omega } (w+t)^2 d A(w+t) \\ 
&=&\int_{h(K)} \lambda_{\Omega -t} (w)^2 d A(w), 
\end{eqnarray*}
where $A$ is the Lebesgue area measure that is also invariant under translations. Due to the uniform convergence of the hyperbolic metric on compacta in Theorem \ref{limdist}, we obtain
\begin{eqnarray*}\label{area2}
\lim_{t \to - \infty} \Ar_h^{\D} (\phi_t(K)) &=& \lim_{t \to - \infty} \int_{h(K)} \lambda_{\Omega -t} (w)^2 dA(w) \\
&=& \int_{h(K)} \lim_{t \to - \infty} \lambda_{\Omega -t} (w)^2 dA(w)  \\
&=& \int_{h(K)} \lambda_{h(\Delta)} (w)^2 dA(w) \\
&=&\Ar_h^{\Delta} (K). 
\end{eqnarray*}

Concerning the hyperbolic $n$-th diameter, suppose that $\phi_t(z_1), ..., \phi_t(z_n)$ is a Fekete $n$-tuple of the compact set $\phi_t(K)$ in $\D$. Then 
\begin{equation}\label{diamdef}
    d_{n,h}^{\D}(\phi_t(K))^{\frac{n(n-1)}{2}}= \prod_{1 \leq \mu < \nu \leq n} \tanh d_{\D}(\phi_t(z_{\mu}), \phi_t(z_{\nu})) .
\end{equation}
From Theorem \ref{monothypdiam}, the limit of the hyperbolic $n$-th diameter, as $t \to - \infty$, exists and we have 
\begin{eqnarray} \label{diam1}
    \lim_{t \to -\infty}  d_{n,h}^{\D}(\phi_t(K))^{\frac{n(n-1)}{2}}&= & \prod_{1 \leq \mu < \nu \leq n}  \lim_{t \to -\infty} \tanh d_{\D}(\phi_t(z_{\mu}), \phi_t(z_{\nu}))  \nonumber\\
    &= &\prod_{1 \leq \mu < \nu \leq n} \tanh d_{\Delta} (z_{\mu}, z_{\nu}) \nonumber\\
    & \leq &d_{n,h}^{\Delta}(K)^{\frac{n(n-1)}{2}},
\end{eqnarray}
where the second equality holds due to Theorem \ref{limdist}. 
Furthermore, suppose now that $y_1, ..., y_n$ is a Fekete $n$-tuple of $K$ in the petal $\Delta$. Then 
\begin{eqnarray}\label{diam2}
d_{n,h}^{\Delta}(K)^{\frac{n(n-1)}{2}} & = &\prod_{1 \leq \mu < \nu \leq n}  \tanh d_{\Delta}(y_{\mu}, y_{\nu}) \nonumber \\ 
&= &\prod_{1 \leq \mu < \nu \leq n} \tanh \left(\lim_{t \to - \infty }   d_{\D}(\phi_t(y_{\mu}), \phi_t(y_{\nu})) \right) \nonumber \\ 
& = &\lim_{t \to -\infty} \prod_{1 \leq \mu < \nu \leq n} \tanh d_{\D}(\phi_t(y_{\mu}), \phi_t(y_{\nu}))  \nonumber \\ 
& \leq & \lim_{t \to - \infty} d_{n,h}^{\D} (\phi_t(K))^{\frac{n(n-1)}{2}}, 
\end{eqnarray}
where the second equality holds due to Theorem \ref{limdist}. 
Combining \eqref{diam1} and \eqref{diam2}, we get $$\lim_{t \to - \infty} d_{n,h}^{\D} (\phi_t(K)) = d_{n,h}^{\Delta}(K). $$
 \qed
 
\section{\bf Condenser Capacity - Proof of Theorem \ref{capacity}}\label{condenser}
Before the proof of Theorem \ref{capacity}, we state the following remark concerning the convergence of the Green function. 

\begin{remark}\label{convergenceGreen}
According to Theorem \ref{limdist} and the connection between the hyperbolic distance and the Green function in a simply connected domain \eqref{greenformula}, for $z, w \in \Delta$
   \begin{eqnarray*}
   \lim_{t \to -\infty} g_{\D}(\phi_t(z), \phi_t(w)) &=& \lim_{t \to -\infty} \log \tanh d_{\D} (\phi_t(z), \phi_t(w))  \\
   &=& \log \tanh d_{\Delta}(z,w) \\
    &=& g_{\Delta}(z,w). 
   \end{eqnarray*}
\end{remark}

As described in Subsection \ref{diameter}, the pairs $(\D, \phi_t(K))$ and $(\Delta, \phi_t(K))$ form condensers in $\mathbb{C}$, for all $t \leq 0$. 
Due to the conformal invariance of condenser capacity \eqref{confinvcapacity} and the fact that $h(\Delta)$ is a horizontal domain, we obtain
   \begin{eqnarray*}
   \Capac (\mathbb{D},\phi_t(K)) &=& \Capac (\Omega,h(K)+t)\\
   &\leq& \Capac(h(\Delta),h(K)+t)\\
   &=&\Capac (h(\Delta),h(K))\\
   &=&\Capac (\Delta,K),
   \end{eqnarray*}
   for all $t\leq 0$. Therefore,
   \begin{equation}\label{cond1}
   \limsup\limits_{t\to -\infty} \Capac (\mathbb{D},\phi_t(K))\le \Capac (\Delta,K).
   \end{equation}
   Now, we need a reverse inequality. We write $\Capac (\mathbb{D},\phi_t(K)) =\Capac (\Omega-t,h(K))$. Since $h(\Delta) \subset \Omega-t$, from the Strong Markov Property for the Green function (Proposition \ref{markovgreen}), for $z, w \in h(K)$, we have
   \begin{equation}\label{markovgreenproof}
       g_{\Omega-t}(z,w)=g_{h(\Delta)}(z,w)+\int\limits_{\partial h(\Delta)}g_{\Omega-t}(\zeta,z)\cdot\omega(w,d\zeta,h(\Delta)).
   \end{equation}
From Remark \ref{convergenceGreen}, when $\zeta \in \partial h (\Delta)$, we have that $g_{\Omega-t}(\zeta,z) \to g_{h(\Delta)} (\zeta,z ) =0 $, as $t \to -\infty$, since the Green function vanishes on the boundary. Hence $g_{\Omega-t}(\cdot, z)$ converges pointwise to $0$ on $\partial h(\Delta)$. For the uniform convergence, we examine the supremum of $g_{\Omega-t}(\cdot, z)$ on $\partial h(\Delta) $, which is attained at the closest point to $z$, in the sense of the hyperbolic geometry in $\Omega-t$. There exists a point $x_t \in  \partial h(\Delta)$ such that 
$$\sup_{\partial h(\Delta)} g_{\Omega-t}(\cdot, z) = g_{\Omega-t} (x_t, z). $$
However, $ \dist (x_t, \partial \Omega-t) \xrightarrow{t \to -\infty} 0$ that implies $\lambda_{\Omega-t}^{\star} (x_t) \xrightarrow{t \to -\infty} + \infty$ for the density of the quasi-hyperbolic metric, and so, the quasi-hyperbolic distance $d_{\Omega-t}^{\star}(x_t, z) \xrightarrow{t \to -\infty} + \infty $. 
From inequality \eqref{relationquasidis} along with the connection between the Green function and the hyperbolic distance \eqref{greenformula}, we obtain 
\begin{equation}\label{ineqgreenquasi}
    \log \tanh \frac{1}{4} d_{\Omega-t}^{\star}(x_t, z) \leq g_{\Omega-t} (x_t, z) \leq  \log \tanh d_{\Omega-t}^{\star}(x_t, z).
\end{equation}
Taking the limit as $t \to -\infty$, we obtain $ g_{\Omega-t} (x_t, z) \to 0. $
Hence $$\lim_{t \to - \infty} \sup_{\partial h(\Delta)} g_{\Omega-t}(\cdot, z)  = 0$$ and we conclude the uniform convergence of $g_{\Omega-t}(\cdot, z) $ on $\partial h(\Delta) $ to $0$. As a result, for every $\epsilon > 0$, there exists a $t_0 \leq 0$, such that for every $t \geq t_0$, $g_{\Omega-t} (\zeta, z) < \epsilon$, for all $\zeta \in \partial h(\Delta)$. 
   Returning back to \eqref{markovgreenproof}, we obtain that for $z, w \in h(K)$, 
   \begin{eqnarray}\label{greenineq}
    g_{\Omega-t}(z,w) & < & g_{h(\Delta)}(z,w)+\int\limits_{\partial h(\Delta)} \epsilon \cdot\omega(w,d\zeta,h(\Delta)) \nonumber\\
    &= &g_{h(\Delta)}(z,w)+ \epsilon \cdot \omega (w, \partial h(\Delta), h (\Delta)) \nonumber\\
    &=& g_{h(\Delta)}(z,w)+ \epsilon.  
    \end{eqnarray}
   The set $h(K)$ is compact, so, it has a unique Green equilibrium measure in each domain. Suppose $\mu_t$ is the Green equilibrium measure of $h(K)$ in $\Omega -t$ and $\mu $ the Green equilibrium measure of $h(K)$ in $h(\Delta)$. Integrating in \eqref{greenineq} with respect to $\mu$, we obtain 
   \begin{eqnarray*}
   \iint_{h(K)^2} g_{h(\Delta)}(z,w) d \mu (z) d\mu (w) + \epsilon \cdot \mu(h(K))^2 &> & \iint_{h(K)^2} g_{\Omega-t}(z,w) d \mu (z) d\mu (w) \\
   &> & \iint_{h(K)^2} g_{\Omega-t}(z,w) d \mu_t (z) d\mu_t (w),
   \end{eqnarray*}
   according to the properties of Green equilibrium measures and bearing in mind that $\mu(h(K))=1$, it follows
   \begin{equation}\label{energies}
      V(h(K), \Omega-t) < V(h(K), h(\Delta)) + \epsilon .
   \end{equation}
   The inequality \eqref{energies} leads to $$\limsup_{t \to -\infty} V(h(K), \Omega-t) \leq V(h(K), h(\Delta))\Rightarrow \limsup_{t \to -\infty} \frac{2 \pi}{\Capac(\Omega-t, h(K))} \leq \frac{2 \pi }{\Capac(h(\Delta), h(K))}.$$
   However $h(K)$ is non-polar, since $K$ is also non-polar, so, $\Capac(\Omega-t, h(K))>0$ and we can conclude that 
   \begin{equation}\label{cond2}
       \liminf_{t \to -\infty} \Capac(\Omega-t, h(K)) \geq \Capac(h(\Delta), h(K)).
   \end{equation}
    Combining \eqref{cond1} with \eqref{cond2} and using the conformal invariance of condenser capacity, we get
    $$ \limsup_{t\to -\infty} \Capac (\mathbb{D},\phi_t(K))\leq \Capac (\Delta,K)\leq \liminf_{t \to -\infty} \Capac(\D, \phi_t(K)).$$
   As a result, $ \Capac(\D, \phi_t(K)) \xrightarrow{t \to -\infty} \Capac (\Delta,K)$.  \qed
   
\begin{remark}
From Theorem \ref{capacity} and \eqref{greencapacity2}, we can conclude the following about the asymptotic behavior of the hyperbolic capacity: 
$$\lim_{t \to - \infty} \caph \phi_t(K) = \caph_{\Delta} K. $$ 

\end{remark}

\section{\bf Non-Regular Backward Orbits - Proof of Theorem \ref{nonreg}}\label{nonregular}

As stated in Subsection \ref{semigroups}, a non-regular backward orbit $\gamma:[0,+\infty)\to\mathbb{D}$ determines a degenerate petal, whose image under the associated Koenigs function of $(\phi_t)$ is a horizontal line. 
If we denote by $\Delta$ this degenerate petal, we get $h(\Delta)=\{h(\gamma(0))+t:t\in\mathbb{R}\}$. 

\begin{figure}[h]
    \centering
    \begin{subfigure}[b]{0.4\textwidth}
    \centering
    \includegraphics[width=\textwidth]{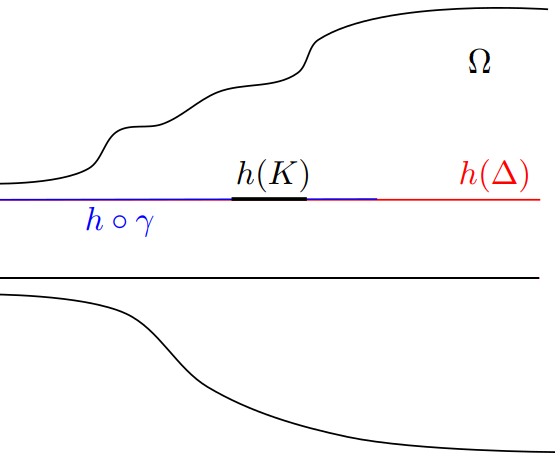}
    \caption{Case (i)}
    \end{subfigure}
    \hfill
     \begin{subfigure}[b]{0.4\textwidth}
    \centering
    \includegraphics[width=\textwidth]{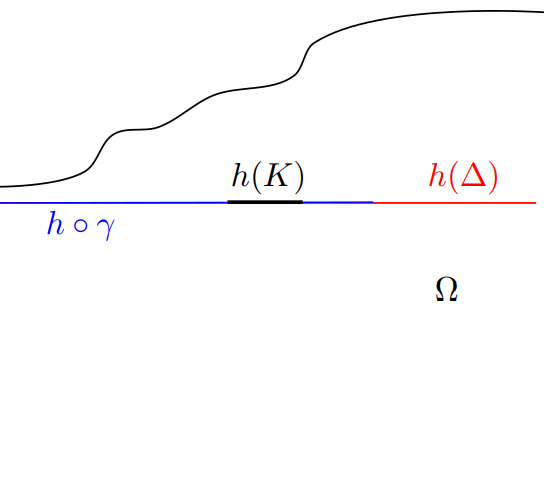}
    \caption{Case (ii)}
    \end{subfigure}
    
    \begin{subfigure}{0.4\textwidth}
    \centering
    \includegraphics[width=\textwidth]{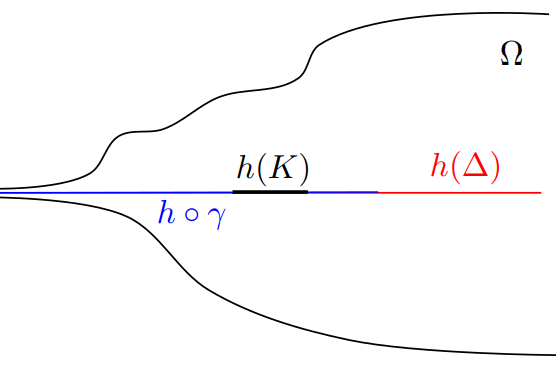}
    \caption{Case (iii)}
    \end{subfigure}
    
    \caption{Non-regular backward orbit of $K$}
    \label{nonregularimage}
\end{figure}

Suppose $K$ is a compact subset of $\Delta$. In order to avoid polar sets, we assume that $K$ is a continuum on $\Delta$. Hence $h(K)$ is a line segment. 
If $K$ has a different form, the proof is analogous. 
In Figure \ref{nonregularimage}, we see all the possible cases when $h(K)$ lies on the image of a degenerate petal. 
 An important property that all the three cases of non-regular backward orbits possess is the fact that
 \begin{equation}
     \lim\limits_{t\to+\infty}\delta_\Omega(h(\gamma(t)))=0,
 \end{equation}
 where $\delta_\Omega(z)=\text{dist}(z,\partial\Omega)$, for any $z\in\Omega$. 

\textbf{Proof of Theorem \ref{nonreg}.} 
Suppose $K$ is a compact non-polar subset of $\Delta$. Then,
\begin{enumerate}[(i)]
    \item Let $z\in\Delta\setminus K$. Then, for $t\le0$, by the conformal invariance of the harmonic measure, we have
    \begin{eqnarray*}
    \omega(\phi_t(z),\phi_t(K),\mathbb{D})&=&\omega(h(z)+t,h(K)+t,\Omega)\\
    &=&\omega(h(z),h(K),\Omega-t).
    \end{eqnarray*}
    It is clear that $h(\Delta)\subseteq\partial\left(\bigcap\limits_{t\le0}(\Omega-t)\right)$ and as a result, $\lim\limits_{t\to-\infty}\delta_{\Omega-t}(h(z))=0,$ while the distance between $h(z)$ and $h(K)$ remains constant, regardless of $t$. 
    Hence, the probability that a Brownian motion starting from $h(z)$ hits $\partial h(K)$ before the boundary of $\Omega -t$, for small values of $t$, is negligible. 
    In conclusion, $\omega(h(z),h(K),\Omega-t)
    \xrightarrow{t\to-\infty} 0.$
    \item Let $z\in\Delta$. By conformal invariance, for $t\le0$, we find that
    $$\lambda_\mathbb{D}(\phi_t(z))=\lambda_\Omega(h(z)+t)=\lambda_{\Omega-t}(h(z))\ge\frac{1}{4\delta_{\Omega-t}(h(z))},$$
    with $\lim\limits_{t\to-\infty}\delta_{\Omega-t}(h(z))=0.$ Consequently, $$\lim\limits_{t\to-\infty}\lambda_\mathbb{D}(\phi_t(z))=+\infty.$$
    \item Let $z,w\in\Delta$ with $z\neq w$. With the help of the conformal invariance of the hyperbolic distance, we have for $t\le0$,
    \begin{eqnarray*}
    d_\mathbb{D}(\phi_t(z),\phi_t(w))&=&d_\Omega(h(z)+t,h(w)+t)\\
    &=&d_{\Omega-t}(h(z),h(w))\\
    &\ge&\frac{1}{4}\log\left(1+\frac{|h(z)-h(w)|}{\min\{\delta_{\Omega-t}(h(z)),\delta_{\Omega-t}(h(w))\}}\right).
    \end{eqnarray*}
    The last inequality clearly implies that $$\lim\limits_{t\to-\infty}d_\mathbb{D}(\phi_t(z),\phi_t(w))=+\infty.$$
    \item Keeping in mind that for any $z,w\in\mathbb{D}$, it is true that 
    $$g_\mathbb{D}(z,w)=\log\tanh d_\mathbb{D}(z,w),$$
    the desired result is a direct corollary of (iii).
    \item Once again, by conformal invariance,
    $$\Ar_{h}^\mathbb{D}(\phi_t(K))=\Ar_{h}^\Omega(h(K)+t)=\Ar_{h}^{\Omega-t}(h(K)),$$
    for all $t\le0$. As we mentioned before, the compact set $h(K)$ can be considered as a line segment. Therefore, trivially, its hyperbolic area with respect to $\Omega-t$ is zero, for all $t\le0$. As a result, we directly get $\lim\limits_{t\to-\infty}\Ar_{h}^\mathbb{D}(\phi_t(K))=0$.
    \item Let $n\in\mathbb{N}$. For $t\le0$, let $\phi_t(z_1),\phi_t(z_2),...,\phi_t(z_n)$ be a Fekete n-tuple for $\phi_t(K)$. Then, 
    $$d_{n,h}^{\D}(\phi_t(K))=\prod_{1\le\mu<\nu\le n} \tanh d_\mathbb{D}(\phi_t(z_\mu),\phi_t(z_\nu))^{\frac{2}{n(n-1)}}.$$ Therefore, using (iii), we get the desired result.
    \item Suppose $\mu_t$ is the equilibrium measure of $h(K)$ in $\Omega-t$ and $\mu $ the equilibrium measure of $h(K)$ in $\Omega$. According to \eqref{greencapacity2},
    \begin{eqnarray*}
        \frac{1}{\Capac (\Omega, h(K)+t)} &=& \iint g_{\Omega-t} (h(z), h(w)) d \mu_t(h(z)) d\mu_t(h(w)) \\
        &\leq& \iint g_{\Omega -t} (h(z), h(w)) d \mu (h(z)) d\mu (h(w)) . 
    \end{eqnarray*}
    The family of domains $(\Omega -t )$ lies in $\Omega$ and according to the Subordination Principle of the Green function \cite[\S\ 4.4]{ransford},
    $$g_{\Omega -t} (h(z), h(w)) \leq g_{\Omega}(h(z),h(w)), $$
    where $g_{\Omega}(h(z),h(w))$ is integrable, as $$\iint g_{\Omega}(h(z),h(w)) d \mu (h(z)) d\mu (h(w)) = V(h(K), \Omega) < + \infty. $$
    As a result, we can apply the reverse inequality in Fatou's Lemma and taking $\limsup$, as $t \to -\infty$, we obtain
    \begin{eqnarray*}
    \limsup_{t \to -\infty}  \frac{1}{\Capac (\Omega, h(K)+t)} &\leq &  \limsup_{t \to -\infty} \iint g_{\Omega -t} (h(z), h(w)) d \mu (h(z)) d\mu (h(w)) \\
        &\leq & \iint \limsup_{t \to -\infty} g_{\Omega -t} (h(z), h(w)) d \mu (h(z)) d\mu (h(w)) \\ 
    &\underset{(iv)}{=}& 0.
    \end{eqnarray*}
     Since the capacity is always non-negative and conformally invariant, we obtain that $$\limsup_{t \to -\infty}  \frac{1}{\Capac (\D, \phi_t(K))} =0 \Rightarrow \liminf_{t \to -\infty} \Capac(\D, \phi_t(K)) = +\infty, $$
    which provides us with the desired result. \qed
\end{enumerate}

\section{Conflicts of interests/Competing interests}

The authors declare that there is no conflict of interest.

\medskip

\bibliographystyle{plain}

\end{document}